\documentclass{amsart}

\usepackage{bm}
\usepackage{graphicx}
\usepackage{amsmath}
\usepackage{hyperref}
\usepackage{wrapfig}
\usepackage{amsfonts}
\usepackage{amsthm}

\numberwithin{equation}{section}
\newtheorem{thm}{Theorem}[section]

\newcommand{\be}{\begin{equation}}
\newcommand{\ee}{\end{equation}}
\newcommand{\la}{\label}

\begin{document}

\title[Efficient algorithms for topological inference on random graphs]{Efficient algorithms for topological inference on random graphs}

\author[Pranav Warman]{Pranav Warman}
\email{rajwarman@hotmail.com}
%
%
\author[Iuliana Teodorescu]{Iuliana Teodorescu}
\email{iuliana@usf.edu}
%
%
\author[Razvan Teodorescu]{Razvan Teodorescu}
\email{razvan@usf.edu}

\address{Department of Mathematics and Statistics, University of South Florida}


\begin{abstract}
In this study, we investigate the problem of classifying, characterizing, and designing efficient algorithms for hard inference problems on planar graphs, in the limit of infinite size. The problem is considered hard if, for a deterministic graph, it belongs to the NP class of computational complexity. A typical example rich in applications is that of connectivity loss in evacuation models for natural hazards management (e.g. coastal floods, hurricanes). Algorithmically, this model reduces to solving a min-cut (or max-flow) problem, with is known to be intractable. The current work covers several generalizations: posing the same problem for non-directed networks subject to  random fluctuations (specifically, random graphs from the  Erd\"os-R\'enyi class); finding efficient convex classifiers for the associated decision problem (deciding whether the graph had become disconnected or not); and the role played by choice of topology (on the space of random graphs) in designing efficient, convex approximation algorithms (in the infinite-size limit of the graph). 
\end{abstract}

\maketitle

\section{Introduction}

In most models of evacuation problems, the basic parameters consist of a given transportation network (e.g. highway system) and  some strategies for incorporating social/behavioral dimensions \cite{1}-\cite{27}. An important problem in need of more in-depth exploration is the sudden loss of connectivity in a basic transportation network due to catastrophic self-driven increase in density, bearing much similarity to the jamming transition in statistical physics. The important distinction between these two phenomena lies in the fact that the jamming (or ``freezing") transition in evacuation models is fully non-equilibrium, brought about by the self-interactions between the transportation agents themselves, and is therefore not a thermodynamic phase transition in a classical sense. In this paper we briefly introduce the problem, identify the main mathematical aspects which allow for efficient algorithms used in predictive models, and outline the analytical results which distinguish ``predictable" random graphs from the others. The discussed is based on  preliminary results describing the case of random Erd\"os-R\'enyi graphs, but the methods outlined can be used for a much larger class of random graphs and problems.  

For the particular example mentioned above, the specific algorithmic problem can be formulated as follows: suppose the random graph which describes the dynamical transportation network is ``close" to becoming disconnected, i.e. it can be represented as two disjoint subsets of vertices connected, on average, by a number $O(1)$ of edges (while the total number of vertices $n \gg 1$). Is there a {\emph{local}} algorithm that can predict with high probability the (topological) phase transition from a connected graph to a disconnected one? 

The locality condition means that the algorithm, run at vertex $v$ from the vertex set of the graph, will require for input only the edge set containing $v$, and will provide output updating only the vertex set $\langle v \rangle$ of nearest neighbors of $v$. As a canonical example (and widespread implementation), algorithms based on random walkers (RW)  are local algorithms (e.g., simple RW, self-avoiding RW, etc.)   

The high probability condition can be restated by stating that the algorithm is a bounded-error random algorithm (the precise definition will be recalled in the next section). This type of algorithms is associated with a formal solution of problems of {\bf{\#P}} hardness (by Stockmeyer's theorem), but only if an {\bf{NP}} oracle is available, and even then the result is not constructive. 

\subsubsection{Graph-theoretic formulation of evacuation models: the example of bipartite approximation}

In the most general situation, a discrete model for dynamic connectivity 
consists of a graph $X(V, E)$ whose vertex set $V$ is fixed (representing 
the sites to be evacuated, the destination sites, and transit nodes), while 
the edge set $E$ may be dynamically allocated (edges may effectively disappear 
due to traffic congestion, or simply as a result of the hazard situation to model). 
Moreover, the model would include variables defined on the vertex set, such as 
occupation numbers (or densities), and corresponding correlation coefficients 
defined on the edge set. The fundamental problem then becomes that of inference 
on such a system of coupled stochastic, dynamic variables, in particular for the case 
of large dimensions ($|V|, |E| \to \infty$). 

As an example, consider the (simpler) case in which the underlying transportation 
network is identified with a bipartite graph, in which the sites to be evacuated (labeled 
``bits" as we explain below) are connected to the destination sites (called ``checks"), but not to each other (and similarly for the check vertices). The nomenclature used originates from the 
field of theoretical coding, in particular Low Density Parity-Check codes (LDPC), where such graphical 
models were used to implement multiple binary checks on sequences of bits.   

To implement the randomness associated with a hazard scenario, we consider an ensemble of such graphs 
with a fixed number of bits $n$, number of checks $m$, bit degree
$q$ and check degree $r$, but selected at random from all the configurations satisfying these conditions.  An from the ensemble is then described
in terms of the $n\times m$ adjacency matrix ${A}_{i}^{\alpha}, i = \overline{1, n}, \alpha = \overline{1, m}$, chosen with uniform probability; thus the ensemble is characterized by the following probability distribution
\begin{eqnarray}
&& p({A})=Z^{-1}\prod_i\delta\left(\sum_\alpha A_i^\alpha - q\right)
\prod_\alpha\delta\left(\sum_i A_i^\alpha - r\right) \label{ens}\\
&& Z=\sum_{{A}}\prod_i\delta\left(\sum_\alpha A_i^\alpha - q\right)
\prod_\alpha\delta\left(\sum_i A_i^\alpha - r\right),\label{ensZ}
\end{eqnarray}
where $Z$ is the total number of elements in the ensemble, and $\delta$ is the Kronecker symbol, $\delta: \mathbb{N} \to \{ 0, 1\}, \delta(0) = 1,  \delta(n)= 0, n \ne 0$.

With these definitions, we can study various expectations (averages) with respect to the {\emph{ensemble}} (\ref{ens})
\begin{eqnarray}
\mathbb{E}_{p} [F({A})] \equiv \sum_{{A}}
p({A}) F({A}).\label{Av}
\end{eqnarray}

Unfortunately, efficient computational implementations of these exact formulas are not available: for example, the (odd cycle transversal) problem (reduction of the original graph to a bipartite one) is in the {\bf{NP}} class. 

This paper is organized as follows: after reviewing the main concepts and results related to graph connectivity detection via local algorithms, we present a constructive approach for the design of efficient random algorithms for the problem of topological inference, specific to given classes of random graphs. We conclude with a discussion of further research topics, such as ergodicity of local random algorithms and its relation to the choice of topology on the space of random graphs (in the infinite-size limit), to be presented in forthcoming publications. 

\section{Complexity classes and efficient approximations}

In this section, we recall some fundamental notions from the theory of computation, along with analytical tools for describing the topological structure of graphs. This short review serves both placing the problem into context, and to establish notational conventions. 

{{We say that a {\emph{decision problem}} $\mathcal{D}$ is in class {\bf{P}} if there exists $k \in \mathbb{N}$ and a deterministic algorithm solving  $\mathcal{D}$ at input size $n$ in computational time $T(n) = O(n^k)$.}}
A {{decision problem}} $\mathcal{D}$ is in class {\bf{NP}} if the problem of verifying any given instance solving $\mathcal{D}$ is in {\bf{P}}.

Equivalently, problems from the class {\bf{NP}} are referred to as being untractable. A canonical example was provided by L. Valiant in 1979, when he showed that every {\bf{NP}}-hard problem can be reduced to the computation of the permanent of  
an arbitrary matrix $A_{n\times n}$, 
$
Per(A) \equiv \sum_{\sigma \in S_n} \prod_{k=1}^n A_{k, \sigma(k)},
$
where  $S_n$ is the symmetric group of $n$ elements.
Another famous {\bf{NP}}-hard problem is that of constraint satisfiability, or 3-SAT: 
Find $\Sigma = \{x_j \}_{j=1}^n$ Boolean variables satisfying 
$
 \wedge_{i=1}^k Q_i(x_a, x_b, x_c), \, x_{a, b, c} \in \Sigma,
$
with $Q$ any logical expression in the 3 variables. 

{{A {\emph{counting problem}} $\mathcal{C}$ is in class {\bf{\#P}} if it is equivalent to counting all instances solving a problem $\mathcal{D} \in $ {\bf{NP}}.
}}
If any other problem $X \in $ {\bf{\#P}} can be reduced by a polynomial-time procedure to $\mathcal{C}$, then $\mathcal{C}$ is a {\bf{\#P}}-complete problem.  

For the class {\bf{\#P}}, a famous result was derived by Stockmeyer, and it involves solvability based on existence of {\emph{oracles}}. 
An {\emph{oracle}} $O$ for class $\mathcal{C}$ is a computational device which can verify if a given instance solves a problem $X \in \mathcal{C}$ (or not), in $O(1)$ computational time. 

In order to describe Stockmeyer's theorem, we also need to introduce the notion of random algorithms: 
A random (or probabilistic) algorithm is a non-deterministic algorithm whose transition from one state to the next (in computational time) is stochastic. 
For random algorithms, the appropriate notion of solvability is based on the concept of statistical sampling with independent events: 
We say that a problem $X$ is in the {\emph{bounded-error}} class $B \mathcal{C}$ if it can be solved by a probabilistic algorithm in class $\mathcal{C}$, whose probability of failure can be made arbitrarily small by repeated application of the algorithm (sampling in statistical sense).  

Then subject to the (strong) requirement of existence of an {\bf{NP}} oracle, the usefulness of random algorithms is revealed in the following theorem:

\begin{thm}(Stockmeyer) 
There is a bounded-error random algorithm that can approximate {\bf{\#P}} using an {\bf{NP}} oracle.
\end{thm}

Nonetheless, in most untractable problems, an {\bf{NP}} oracle is not available, so it is necessary to consider extensions of this theorem to the case of weaker conditions, and perhaps weaker guarantees for an approximate solution.

\subsection{Topological characteristics of graphs and associated generating functions} For a generic graph $G(V,E)$, the Ihara (or vertex) zeta function of $G$ is defined as the complex-valued function with respect to a formal variable $u \in \mathbb{C}$,
\be \la{ihara}
\zeta_V(u,G) = \prod_{\gamma \in \mathcal{P}(G)} \left [ 1-u^{l(\gamma)}\right ]^{-1}, \quad u \in \mathbb{C},   
\ee
where $\mathcal{P}(G)$ is the set of all ``prime" loops in the graph (no tail, no backtracking, and not a power of another loop).  An immediate generalization follows from the notion of ``graph inflation": if
$w : E(G) \to \mathbb{N}$ is an integer weight function ($w(e)$ is a positive integer associated with 
 edge $e$), then its weighted Ihara zeta function is
\be \la{ihara}
\zeta_V(u,G,w) = \prod_{\gamma \in \mathcal{P}(G)} \left [ 1-u^{l(\gamma, w)}\right ]^{-1}, \quad 
l(\gamma, w) = \sum_{e_i \in \gamma} w(e_i),
\ee
with $l(\gamma, w)$ the weighted length of loop $\gamma$. We have the following theorem:
\be
\zeta_V(u,G,w) = \zeta_V(u, G_w, 1),
\ee
where $G_w$ is the ({\emph{inflated}}) graph obtained from $G$ by replacing edge $e$ of weight $w(e)$ with 
$w(e) $ edges of weight one, connected by $w(e) -1$ new degree two vertices. 

The graph seta functions are useful as generating functions for topological characteristics of the graph such as:
$R_G \ge 0$ is the {\emph{radius of convergence}} of the 
graph, defined as the largest positive real such that $\zeta_V(u,x) < \infty \,\, \forall \,\, |u| < R_G$;
 $N_m$, the number of closed loops (including multiples of simple loops) of length $m$ in $G$; and $\kappa_G$, the number of spanning trees in $G$. 

An exact analytical expression of the Ihara zeta functions produces the topological characteristics by the following formulas: 
$$
N_m  = \frac{1}{m!} \left [ \frac{\partial }{\partial u} \right ]^m_{u = 0} \frac{\partial \log \zeta }{\partial \log u}, 
$$
$$
\left [ \frac{\partial }{\partial u} \right ]^{n}_{u = 1} \zeta^{-1}_V(u,G) = n! (-1)^{n+1}2^n(n-1)\kappa_G,
$$
Moreover, there is an asymptotic expansion formula as $m \to \infty$, 
$$\pi(m) = \Delta_G \frac{R^{-m}_G}{m} + O(m^{-2}),$$
where $\pi(m)$ is the number of prime loops of length $m$ in $G$, and $\Delta_G$ is the greatest common denominator of lengths of prime loops in $G$. 

Therefore, it may appear that a good way to predict average topological characteristic of random graphs consists of computing the averaged zeta function, and then compute directly the quantity of interest. However, the computational complexity for the functional evaluation of zeta functions is  also in the {\bf{NP}} class, which means that this approach becomes prohibitively complex in the limit of large graphs (precisely the limit of interest). Furthermore, there is no known algorithm to compute $\zeta_V$ by local algorithms. A tantalizing result is provided by the Ihara-Hashimoto-Bass theorem, expressing the zeta function by means of a determinant representations: 
\begin{thm}
For a graph $G(V, E)$ and Euler characteristic $\chi(G) = |E| - |V|$, 
$$
\zeta^{-1}_V(u,G) = (1-u^2)^{-\chi(G)} {\sf{Det }} [I - uA + u^2Q],
$$
where $Q$ represents the degree matrix.  
\end{thm}

Any determinantal representation is of great interest from a computational point of view, given that computation of the determinant can be performed in polynomial time $O(N^3)$. There are, however, two significant difficulties to consider: the first is that the Ihara-Hashimmoto-Bass theorem gives a functional relation with a generic parameter $u \in \mathbb{C}$; the second is finding a local algorithm for computation of the determinant of the $|V| \times |V|$ matrix $I - u A + u^2 Q$. While no algorithm is known to solve both of these problems, the subproblem calling for local stochastic algorithms that effectively compute functional determinants leads us to consider diffusion-like stochastic processes such as Random Walks. 

\subsection{Analysis of Random Walks on graphs}

Given a graph $G(V, E)$ with discrete Laplacian $\Delta$, the evolution of random walkers in time $t \in \mathbb{N}$ is given by 
%
%
\be \la{discrete}
\frac{\Delta p_t(v)}{\Delta t} \equiv p_{t+1}(v) - p_t(v) = \Delta p_t(v).
\ee
%
%
Introducing the two-point correlation function for the random walker, defined by the conditional probability 
\be
C_2(v, w; t, s) \equiv \mathbb{P}( v, t | w, s), 
\ee
we have similarly the formal solutions 
\be
p_n (v) = \sum_{ w \in V} C_2(v, w; n, 0) p_0(w) = \sum_{w \in V} \left [ (\mathbb{I} + \Delta)^n \right ]_{v, w} p_0(w) = 
[ (A_G^n) p_0]_v,  
\ee
where $p_0$ is the vector of probabilities at time $t = 0$,  and  $A_G$ is a scaled adjacency matrix of $G$. Topological characteristics of the graph, such as the number of loops of length $n$, containing the vertex $v$, are expressed by the diagonal elements of matrices $A_G^n$. 

Recall that the continuum case corresponds to the limit of a regular graph with even degree $d=2n$, in the limit of edge length going to zero, and is given (in free space $\mathbb{R}^{d/2}$) by the Gaussian kernel 
\be
C_2(\vec r_1, \vec r_2;  t, s)  = K_{t, s}(\vec r_1, \vec r_2) = \frac{1}{\sqrt{(2\pi)^{d/2}\cdot 2(t-s)}} e^{-\frac{||\vec r_1 - \vec r_2 ||^2}{2(t-s)}}, \quad t > s. 
\ee 
The ``diagonal" correlation function will then produce the prefactor (ignoring the dependence on the diffusion time variable) $C_2(\vec r, \vec r; t, s) \sim [(2 \pi)^{d/2}]^{-\frac{1}{2}}$, which is precisely the volume of the elementary cell of the reciprocal graph for $G$, in the continuum limit. 

For a proper discrete example, we consider the one-dimensional case of an asymmetric RW. For this particular case, we take $0\le p, q \le 1, p+ q = 1$, and denote the evolution operator by $T = pW + qW^*$, where $W$ represents a ``step to the right" with probability $p$, and its adjoint $W^*$ represents  a ``step to the left" with probability $q$, on the lattice $\mathbb{Z}$. Therefore choosing for initial distribution at ``time" $n = 0$ to be the pure state $\{p_n(0) = \delta_{n, 0}\}_{n \in \mathbb{Z}}$, 
we conclude that Trace$(T^n)$ computes the probability $P_n$ for the random walker starting from $0$ at $n = 0$ to be found at $0$ at some time $n > 0$, which means that the function, 
$$
G_T(z) \equiv  \sum_{n = 0}^{\infty} \frac{P_n}{z^{n+1}}
$$
is a generating function for the probabilities of return in $n$ steps, $\{ P_n \}_{n \ge 0}$. To evaluate this generating function, remark that the spaces $\ell^2(\mathbb{Z})$ and $L^2(\mathbb{T})$ are isomorphic, so that we can map $W \to w = e^{i \theta} \in \mathbb{T}$ (the unit circle), 
so that (for $|z|$ sufficiently large),
$$
G_T(z) = \frac{1}{2\pi i} \oint_{\mathbb{T}} \frac{1}{w[z - (pw+qw^{-1})]} dw = \frac{1}{\sqrt{z^2 - 4 pq}}
$$
The quantity analogous to $C_2(0, 0; n, 0)$ is then obtained by computing the integral 
$$
P_n = \frac{1}{2\pi i}\oint_{\mathbb{T}} z^n G_T(z) dz. 
$$

%
%
%

While using independent random walkers can produce characteristic quantities (describing, for instance, the reciprocal graph of a $d-$regular lattice), it is not sufficient for designing efficient algorithms, basically because of the independent R.W. being ``memoryless". A  more efficient procedure would require R.W. to interact (either with each other, or each with its own path). Mutually-avoiding random walkers and self-avoiding random walkers are two such examples. Computationally, however, the evaluation of correlation functions for interacting diffusions becomes a much more involved problem, for which exact results are not available. We propose to circumvent this difficulty by considering a limiting case for the estimation of probability functions for interacting diffusions in the large deviations limit. To that end, we review some of the basic elements of the theory, in the next section.  


\subsection{Large Deviations Principle for i.i.d.r.v.}
We recall that 
the Cram\'er functional (or ``rate function") for large deviations of the sample mean for a random variable $X$ is defined as
$$
I(x) := \max_{t > 0} \, \left [t x - \ln(m_X(t)) \right ],
$$
where $m_X(t)$ is the moment-generating function of $X$, 
$$
m_X(t) = \mathbb{E}(e^{tX}).
$$

The Large Deviations Principle  states that the probability of ``large deviations"
$$
P\left (\frac{1}{n}\sum_{k=1}^n X_k \ge x \right ) \sim e^{-n I(x)}, 
$$
where $X_1, X_2, \ldots, X_n$ is a sample of  i.i.d. r.v., and $I(x)$ is the rate function defined earlier, and $``\sim"$ means that the probability
is determined only up to an overall normalization factor.
For the case of interacting diffusions, assume that the sequence of random variables 
$Z_n, n = 1, 2, \ldots$ from the space $\Sigma$ have distributions 
$d P_n$ and moment-generating functions $m_n(t) = \mathbb{E}_{P_n}(e^{tZ_n})$.  
If the limit
\be \la{fenergy}
\Lambda(t) = \lim_{n \to \infty} n^{-1}\log m_n(t)
\ee
exists, is convex and bounded from below, and that 
\be \la{rate}
Q(x) = \sup_{t}[xt - \Lambda(t)]
\ee
is a well-defined {\emph{rate function}} (bounded from above,
lower semi-continuous and has compact level sets), then 
$\{ Z_n \}$ satisfies the Large Deviations Principle with rate  $Q$:
\be \la{one}
\begin{array}{c}
\lim \sup  \log P(Z_n \in C ) \le -n Q(C), \\ 
\lim \inf \log P(Z_n \in O ) \ge -n Q(O), \,\!\, 
\end{array}
\ee
where $n \to \infty$, the sets $C, O$ are closed and open, respectively, 
and $Q(S)$ is by definition 
\be
Q(S) \equiv \inf_{x \in S} Q(x), \,\, (\forall) \,\, S \subset \Sigma.
\ee
This is known as the G\"artner-Ellis theorem \cite{G-E}. 
For many situations, (\ref{one}) imply that for any set $S$, 
\be \la{ge}
n^{-1} \log P(Z_n \in S) \to - Q(S).
\ee
 
The  G\"artner-Ellis theorem provides the theoretical tools necessary for studying 
the long-time behavior of interacting systems, when their equilibrium limit 
without interactions is known. For example, consider the case of $N$ 
interacting diffusions, described by the coupled stochastic equations
\be \la{ito}
d \eta_i = \left [\sum_{j \ne i} f_i(\eta_i-\eta_j) + g(\eta_i) \right ] dt + d\xi_i,
\ee
where the first terms corresponds to inter-particle interactions, 
the second is a single-particle forcing, and the last  represents
uncorrelated Brownian processes. 

We can use (\ref{ge}) to derive the rate function 
for the  $N$ interacting diffusions (\ref{ito}), in the following way:  
denoting by 
$F_i =\sum_{i \ne j} f(\eta_i, \eta_j) +  g(\eta_i)$, and by $P_N, P_N^{(0)}$ the 
distribution laws for the interacting diffusions $\eta_i$ and non-interacting 
$\xi_i$, respectively, we obtain 
\be
d P_N = 
e^{ \left [ 
-\frac{1}{2}\sum_{i=1}^N \int_0^T (F_i)^2 dt -
\sum_{i=1}^N \int_0^T F_i  d \eta_i(t)
\right ]}
d P_N^{(0)}.
\ee
This formula can be used to approximate the joint distribution of interacting diffusions
starting from their non-interacting (independent) joint distribution. 

%

\subsection{Ergodicity as a condition for numerical convergence of efficient algorithms}

The remaining theoretical result that is required in order to ensure the efficiency of the numerical simulation procedure is known as the von Neumann (or $L^2$-) ergodic theorem, formulated below:

\begin{thm}
Consider the probability space $(\Omega, \Sigma, \mathbb{P})$, where $\Omega \subset \mathbb{R}$, $\Sigma$ is the Borel algebra of  $\Omega$, and $\mathbb{P}$ is a proper distribution on $\Omega$, and the semi-flow $T_n: \mathbb{N} \times \Omega \to \Omega$ is ergodic, i.e. $T_n = T^n$, 
\begin{itemize}
\item[(i)] $\mathbb{P}(A) = \mathbb{P}(T_n(A)), \,\, \forall A \in \Sigma$
\item[(i)] $T_n^{-1}(A) \in \Sigma \,\, \forall A \in \Sigma$
\item[(i)] $T_n(A) = A$ only if $\mathbb{P}(A) = 1$ or if $\mathbb{P}(A) = 0$.
\end{itemize}  
If $U_n$ denotes the operator defined by $[U_n (h)](\theta) = h(T_n(\theta))$, for all functions $h$ square-integrable on $\Omega$, then the series of approximations 
$$
h_n(x) \equiv \frac{1}{n} \sum_{k=1}^n U_n(h(x)) 
$$
converges to the ensemble average of $h(x)$, in $L^2-$norm, with error less than $O(n^{-2})$. 
\end{thm}
 
For our purposes, this theorem provides a criterion for the algorithm to converge (as a process in the {\emph{computational time}} $n \in \mathbb{N}$) to the unique limit which is the ensemble average. Another (perhaps more interesting, and definitely much more challenging) problem to be considered is the dependence of convergence behavior on the choice of topology (and therefore, Borel algebra) on the space of graphs. While for finite-size graphs the limiting behavior does not depend on the choice (although the rate of convergence might), for ensembles of graphs in the infinite-size limit the situation is much more complicated. We will address this interesting aspect in a future publication. 

\section{Designing efficient approximation algorithms: the method of reduction to critical statistical models}
 
As noted by several authors in recent studies on performance of ``deep learning" algorithms \cite{new}, successful machine learning (as dynamical processes in the learning time) explores the data structures at many diverse characteristic ``scales", which allows for fast convergence on the ``critical" subset. The analogy with the variational renormalization group is much deeper than it may seem, because the very existence of a critical data subset implies the existence of a critical manifold, and therefore universal scaling behavior. From this point of view, it is possible to obtain a different classification of convergent algorithms (for given classes of problems), in the same way as critical statistical mechanical models are uniquely described by the corresponding renormalizable minimal model (e.g., the ``$\phi^4$" theory, tri-critical point, edge singularity, etc.). 

We propose to take this observation one step further and explore the possibility of designing an efficient approximating algorithms by coupling a given class of random graphs to its corresponding minimal theory, in the sense afforded by the Knizhnik-Polyakov-Zamolodchikov relations, for critical statistical models in the presence of fluctuating geometry. The role of fluctuating geometry is played by the random graph (in the $V$, $E$ $\to \infty$ limit), and the minimal model which is to be associated with a given type of graph is available from the literature. Moreover, given that we are only interested in topological properties (and transitions), the corresponding class of critical models is wider than for the case of metric-dependent problems. 

In essence, what this procedure provides is the correct choice of stochastic process to be implemented on the given graph (generalizing the simple random walk discussed earlier), and which will produce graph characteristics corresponding to relevant perturbations around the critical manifold. 

As a readily-available example well covered by the literature, we can analyze the problem of ``perfect matchings" on a regular lattice, and its known solution provided by a minimal theory with irreducible representations $j=1/2$ for the group $SU(2)$ (also known as fermionic theory). The stochastic process implemented by free fermions leads to the effective computation of the partition function of the theory (a Pfaffian, in general), and equivalent criterion for perfect matchings given by the Fisher, Kasteleyn, and Temperley (or FKT) algorithm. 

A comprehensive list of minimal critical models associated with certain classes of random graphs (planar or otherwise) is the subject of a forthcoming publication, while run-time experiments for implementation of this procedure will be made available on  GitHub by the  authors, as the implementation codes become available. 

\section{Acknowledgements}

The authors wish to thank the USF Center for Complex Data Systems for support in the computational phase of this project.

\end{document}